\documentclass[10pt]{amsart}

\usepackage{amsmath,amssymb,amsthm,mathtools}
\usepackage{fullpage}
\usepackage{hyperref}
\hypersetup{hidelinks}
\usepackage{enumitem}

\newtheorem{theorem}{Theorem}[section]
\newtheorem{corollary}[theorem]{Corollary}
\newtheorem{proposition}[theorem]{Proposition}

\theoremstyle{remark}
\newtheorem{remark}[theorem]{Remark}
\newtheorem{example}[theorem]{Example}

\DeclareMathOperator{\Irr}{Irr}
\DeclareMathOperator{\End}{End}
\DeclareMathOperator{\Hom}{Hom}
\DeclareMathOperator{\Tr}{Tr}
\DeclareMathOperator{\Span}{span}
\DeclareMathOperator{\LocMax}{LocMax}
\DeclareMathOperator{\Res}{Res}

\newcommand{\C}{\mathbb C}
\newcommand{\PP}{\mathbb P}
\newcommand{\norm}[1]{\lVert #1\rVert}
\newcommand{\ip}[2]{\langle #1,#2\rangle}

\usepackage{xcolor}

\title{Incidence Geometry of Local Maxima:\\ Howe Duality and Branching}

\author{Pin-Chi Hung}
\address{Department of Mathematics, Soochow University, Taipei City, Taiwan}
\email{pinchihung1111@gmail.com}

\author{Ming-Hsuan Kang}
\address{Department of Applied Mathematics, National Yang Ming Chiao Tung University, Hsinchu, Taiwan}
\email{kmsming@gmail.com}

\date{}

\begin{document}

\begin{abstract}
For a finite-dimensional unitary representation of a finite or compact group, consider the mean squared matrix coefficient
on the unit sphere. Its local maxima are the vectors generating irreducible cyclic subrepresentations.
We show that for compatible group actions this variational picture retains the incidence structure of
representation theory. For two commuting actions, intersections of projectivized local-maximal
components recover the weighted support graph of the product representation, including its multiplicity
spaces; Howe duality is the multiplicity-one matching case. More generally, for pairwise commuting
actions, and for subgroup chains, every positive linear combination of the corresponding functionals has
exactly their common local maxima, while every other critical point has a positive Hessian direction.
On a fixed irreducible module at the top of a subgroup chain, the simultaneous projective
local-maximal components are indexed by vertex paths in the weighted Bratteli diagram. Each component
is the product of the projective space of the bottom representation and the projectivized branching
multiplicity spaces. In the multiplicity-free case with one-dimensional bottom
type, these reduce to the Gelfand--Tsetlin lines. Peter--Weyl, Frobenius reciprocity, Young branching and
Schur--Weyl duality provide concrete instances of the general picture.
\end{abstract}

\maketitle


\section{Introduction}\label{sec:intro}

Several standard constructions in representation theory amount to making compatible choices of
irreducible pieces. Howe duality matches irreducibles for two commuting actions; branching records how
an irreducible decomposes along a subgroup chain; and Gelfand--Tsetlin bases arise when successive choices
are multiplicity-free. We ask whether these discrete structures can be recovered as the extremal
geometry of one natural scalar functional.

For a finite or compact group $K$ acting unitarily on a finite-dimensional complex inner-product space $\Omega$, let
\[
 \Phi_K(v)=\int_K |\langle v,kv\rangle|^2\,d\mu_K(k),
 \qquad \|v\|=1,
\]
where $\mu_K$ is normalized Haar measure (the normalized counting measure when $K$ is finite).
All actions of compact groups are assumed continuous. In
\cite[Theorem~3 and Corollary~4]{HungKang2026}, we proved that the local maxima of $\Phi_K$ are exactly
the vectors whose cyclic span is irreducible, and that their connected components are indexed by
the irreducible types occurring in the representation. The compact-group extension is recorded in
\cite[Remark~2]{HungKang2026}. Each component consists of the unit vectors whose cyclic span is
isomorphic to the corresponding irreducible representation.
For a single subgroup inclusion, \cite[Theorem~6]{HungKang2026} also identifies intersections of
projective local-maximal components with products of projective spaces recovering branching
multiplicities. The present paper extends this incidence picture to commuting actions and iterated
subgroup chains, and proves a positive-sum principle for these compatible conditions.

First let $G$ and $H$ be finite or compact groups acting unitarily on $\Omega$, with commuting actions.
The product decomposition has the form
\[
 \Omega\cong\bigoplus_{\pi,\sigma}
 V_\pi\otimes W_\sigma\otimes M_{\pi,\sigma}.
\]
Here $\pi$ and $\sigma$ range over the irreducible
representations of $G$ and $H$, respectively;
$V_\pi$ and $W_\sigma$ are their representation spaces,
and $M_{\pi,\sigma}$ is the multiplicity space of
$V_\pi \otimes W_\sigma$ in $\Omega$.
Representation theory records this by a weighted bipartite support graph. Its left and right vertices are the irreducible
$G$-types and $H$-types occurring in $\Omega$, respectively. An edge $(\pi,\sigma)$ is
present when $M_{\pi,\sigma}\ne0$, with weight $\dim M_{\pi,\sigma}$.
We show that this graph is visible
directly in the two local-maximal decompositions. After quotienting by scalar phase, the $\pi$-component
for $G$ meets the $\sigma$-component for $H$ exactly when the edge occurs, and the intersection is
\[
 \PP(V_\pi)\times\PP(W_\sigma)\times\PP(M_{\pi,\sigma}).
\]
Thus incidence recovers not only support and multiplicity but also the projective family of embedded
copies on each edge (Theorem~\ref{thm:incidence}).
In the finite-dimensional setting considered here, the Howe duality condition
\cite{Howe1989a,Howe1989b,GoodmanWallach} means that the algebras generated by the two actions are
each other's commutants. This is precisely the multiplicity-one matching case. Then $\Omega$ is a sum
of paired blocks $V_\pi\otimes W_\sigma$, and on each block the two functionals coincide after
rescaling by $\dim V_\pi$ and $\dim W_\sigma$. Their common local maxima are exactly the unit
decomposable tensors $x\otimes y$ in the individual paired blocks
(Corollary~\ref{cor:howe}).

Our second result is a variational compatibility theorem. Let $K_1,\ldots,K_s$ be finite or
compact groups acting unitarily on $\Omega$. Assume either that their actions commute pairwise,
or that the groups form a chain $K_1\le\cdots\le K_s$ and their actions are restrictions of a
single unitary representation of $K_s$ on $\Omega$. Then, for arbitrary positive
weights $a_1,\ldots,a_s$,
\begin{equation}\label{eq:intro-positive-sum}
 \LocMax\!\left(\sum_{i=1}^s a_i\Phi_{K_i}\right)
 =\bigcap_{i=1}^s\LocMax(\Phi_{K_i}).
\end{equation}
Here $\LocMax(F)$ denotes the set of local maxima of a function $F$ on the unit sphere of $\Omega$.
Moreover, at every critical point of the sum that is not a local maximum, the Hessian
on the unit sphere has a strictly positive direction
(Theorem~\ref{thm:positive-sum}). Thus compatible irreducibility conditions can be imposed by a single
quartic functional without creating spurious local maxima. This compatibility is essential:
the equality can fail even for two order-two groups
acting on $\mathbb{C}^2$;
see Example~\ref{ex:necessity-of-compatibility}.

For a chain $G_0\leq\cdots\leq G_n$ of finite or compact groups, fix an irreducible type
$\lambda_n\in\Irr(G_n)$ with a finite-dimensional unitary representation space $V_{\lambda_n}$. For
$0\le j<n$ we likewise write $V_{\lambda_j}$ for a representation space of $\lambda_j\in\Irr(G_j)$.
Restriction along the chain is recorded by the weighted Bratteli diagram:
its vertices at level $j$ are the elements of
$\Irr(G_j)$, and there is an edge $\lambda_j\to\lambda_{j-1}$ when $V_{\lambda_{j-1}}$ occurs in the restriction of
$V_{\lambda_j}$ to $G_{j-1}$, with weight equal to its multiplicity. A path below $\lambda_n$ is a sequence
$\lambda_n\to\lambda_{n-1}\to\cdots\to\lambda_0$ of such edges.
Consider the lines $[v]\in\PP(V_{\lambda_n})$ whose unit
representatives are local maxima of every $\Phi_{G_j}$, $0\le j\le n$, on the unit sphere of
$V_{\lambda_n}$. We
show that the connected components of this simultaneous projective local-maximal locus are indexed by the
paths below $\lambda_n$, and that the component of the path $\lambda_n\to\cdots\to\lambda_0$ is naturally
isomorphic to
\[
 \PP(V_{\lambda_0})\times\prod_{j=1}^{n}
 \PP\bigl(\Hom_{G_{j-1}}(V_{\lambda_{j-1}},\Res^{G_j}_{G_{j-1}}V_{\lambda_j})\bigr)
\]
(Theorem~\ref{thm:bratteli}). Each projectivized Hom space parametrizes the embedded irreducible copies
available at one restriction step, while $\PP(V_{\lambda_0})$ records the final line.

By the positive-sum principle, for any $a_0,\ldots,a_n>0$ these are also exactly the
projective local-maximal components of $\sum_{j=0}^{n}a_j\Phi_{G_j}$.
When every edge below $\lambda_n$ has weight one and every
occurring bottom type is one-dimensional, the components
reduce to points and give the lines of the Gelfand--Tsetlin basis
of $V_{\lambda_n}$. Choosing unit representatives gives an orthonormal basis, unique up to
independent phases, whose vectors generate an irreducible $G_j$-submodule for every $j$.
Thus the pair and chain theorems describe the edge and path
forms of the same principle: local maximality selects an
irreducible factor and leaves its multiplicity space as
the residual variable.

The incidence geometry and the connected components of local-maximal loci realize
support and branching data, including the projective families arising from multiplicities.
Related purity extrema occur in generalized entanglement, where Lie-algebraic purity detects generalized coherent states
\cite{BarnumKnillOrtizViola2003}; the group-averaged functional used here is different and retains
isotypic multiplicity and branching data. Section~\ref{sec:nested} gives the residual-module
factorization. Sections~\ref{sec:incidence}--\ref{sec:bratteli} establish the structural results, and
Section~\ref{sec:examples} treats Peter--Weyl, Frobenius reciprocity, Young branching and Schur--Weyl
duality.

\section{The residual multiplicity module}\label{sec:nested}

\paragraph{Notation.}
For a finite group $K$ we write $\Irr(K)$ for its set of irreducible
representations and $V_\alpha$ for a representation space of $\alpha\in\Irr(K)$, of dimension $d_\alpha$.
Let $K$ act unitarily on a finite-dimensional complex inner-product space $E$, with inner product linear in
the second variable, and let $S(E)$ be its unit sphere. Define
\[
 \Phi_K(v)=\frac1{|K|}\sum_{k\in K}|\ip{v}{kv}|^2,\qquad v\in S(E),
\]
and write $\LocMax(\Phi_K)\subset S(E)$ for the set of its local maxima on $S(E)$. For an operator $A$ on
$E$ put
\[
 \mathcal R_K(A)=\frac1{|K|}\sum_{k\in K}kAk^*,
 \qquad
 \norm{A}_{\mathrm{HS}}^2=\Tr(A^*A).
\]
For $P_v=vv^*$, the rank-one projection onto $\C v$, one has
\begin{equation}\label{eq:avg-gradient}
 \Phi_K(v)=\norm{\mathcal R_K(P_v)}_{\mathrm{HS}}^2,
 \qquad
 T_K(v)=\frac{\mathcal R_K(P_v)v}{\norm{\mathcal R_K(P_v)v}}.
\end{equation}
Here gradients are taken with respect to the real inner
product $\operatorname{Re}\langle\cdot,\cdot\rangle$ on $E$.
By \cite[Lemma~1(d)]{HungKang2026},
the ambient gradient of the homogeneous quartic extension
of $\Phi_K$ to $E$ is $4\mathcal R_K(P_v)v$, so $T_K$ is its
normalized ambient-gradient map.

\paragraph{Results used from \cite{HungKang2026}.}
Write the $\alpha$-isotypic component of $E$ as
$V_\alpha\otimes N_\alpha$, with $N_\alpha=\Hom_K(V_\alpha,E)$ carrying the inner product
$d_\alpha^{-1}\Tr(\varphi^*\psi)$, so that the evaluation map $V_\alpha\otimes N_\alpha\to E$ is an
isometry, and regard the $\alpha$-component of $v$ as
an element of $V_\alpha\otimes N_\alpha$, that is, as a $d_\alpha\times\dim N_\alpha$ matrix. If its
singular values are $\sigma_{\alpha,j}$ and $p_{\alpha,j}=\sigma_{\alpha,j}^2$ (the \emph{masses} of $v$, which
sum to $1$), then \cite[Lemma~1(c)]{HungKang2026} gives
\begin{equation}\label{eq:singular-formula}
 \Phi_K(v)=\sum_{\alpha,j}\frac{p_{\alpha,j}^2}{d_\alpha}.
\end{equation}
By \cite[Theorem~3]{HungKang2026},
\begin{equation}\label{eq:criterion}
 v\in\LocMax(\Phi_K)
 \quad\Longleftrightarrow\quad
 \Span\{kv:k\in K\}\text{ is irreducible}.
\end{equation}
Moreover, by \cite[Corollary~4]{HungKang2026}, $\LocMax(\Phi_K)$ is the disjoint union,
over the $\alpha$ occurring in $E$, of the connected sets
$\mathcal C_\alpha^K$ of unit decomposable tensors $x\otimes n$ with $x\in V_\alpha$ and $n\in N_\alpha$; on
$\mathcal C_\alpha^K$ one has $\Phi_K=1/d_\alpha$. We call $\mathcal C_\alpha^K$ the
\emph{$\alpha$-component} of $\LocMax(\Phi_K)$.
For $\alpha\in\Irr(K)$ not occurring in $E$, we set $\mathcal C_\alpha^K=\varnothing$.

\paragraph{Commuting actions.}
Now let $G$ and $H$ act unitarily and commute on $\Omega$. For
$\pi\in\Irr(G)$ let $\Omega[\pi]$ be the $\pi$-isotypic component and write
\begin{equation}\label{eq:theta-module}
 \Omega[\pi]\cong V_\pi\otimes\Theta(\pi),
 \qquad
 \Theta(\pi)=\Hom_G(V_\pi,\Omega).
\end{equation}
We equip $\Theta(\pi)$ with the inner product
$\ip{\varphi}{\psi}_{\Theta(\pi)}=d_\pi^{-1}\Tr(\varphi^*\psi)$, where $d_\pi=\dim V_\pi$; for this
choice the evaluation map $V_\pi\otimes\Theta(\pi)\to\Omega[\pi]$, $x\otimes\varphi\mapsto\varphi(x)$,
is unitary. Since $H$ commutes with $G$, it acts on $\Omega[\pi]$ as $I\otimes\eta_\pi$ for a unitary
representation $\eta_\pi$ of $H$ on $\Theta(\pi)$. We write $\Phi_{\Theta(\pi)}$ and $T_{\Theta(\pi)}$ for the
functional and the normalized ambient-gradient map of the $H$-module $(\Theta(\pi),\eta_\pi)$.

\begin{proposition}[Nested variational decomposition]\label{prop:nested}
Let $\mathcal C_\pi^G$ be the $\pi$-component of $\LocMax(\Phi_G)$. For unit vectors
$x\in V_\pi$ and $m\in\Theta(\pi)$,
\begin{align}
 \Phi_H(x\otimes m)&=\Phi_{\Theta(\pi)}(m),\label{eq:nested-value}\\
 T_H(x\otimes m)&=x\otimes T_{\Theta(\pi)}(m).\label{eq:nested-map}
\end{align}
Consequently the local maxima of $\Phi_H$ on $\mathcal C_\pi^G$ are exactly the vectors $x\otimes m$ for which $m$
generates an irreducible $H$-subrepresentation of $\Theta(\pi)$.
\end{proposition}

\begin{proof}
Since $H$ acts on $V_\pi\otimes\Theta(\pi)$ by $I\otimes\eta_\pi(h)$,
\[
 \ip{x\otimes m}{(1,h)(x\otimes m)}=\ip{m}{\eta_\pi(h)m}.
\]
This proves \eqref{eq:nested-value}. Similarly,
\[
 \mathcal R_H(P_{x\otimes m})
 =P_x\otimes\mathcal R_{\Theta(\pi)}(P_m),
\]
and applying both sides to $x\otimes m$
and normalizing gives \eqref{eq:nested-map}.

Since $\Phi_H(x\otimes m)=\Phi_{\Theta(\pi)}(m)$
is independent of the unit vector $x$,
the vector $x\otimes m$ is a local maximum of
$\Phi_H$ on $\mathcal C_\pi^G$ if and only if $m$ is a local
maximum of $\Phi_{\Theta(\pi)}$ on $S(\Theta(\pi))$.

The last assertion follows from \eqref{eq:criterion}
applied to the $H$-module $\Theta(\pi)$.
\end{proof}

\begin{remark}[Compact groups]\label{rem:compact}
The averaging, singular-value, and local-maximum results recalled above extend to compact
groups acting continuously and unitarily on finite-dimensional spaces, as recorded in
\cite[Remark~2]{HungKang2026}. Replacing finite averages by integrals against normalized Haar
measure leaves Schur orthogonality and the isotypic matrix model unchanged. The arguments below
therefore apply to compact groups as well. We use this in Section~\ref{sec:SW} for $U(n)$.
\end{remark}

\section{Incidence geometry}\label{sec:incidence}

For $\sigma\in\Irr(H)$ write $W_\sigma$ for its representation space and $e_\sigma=\dim W_\sigma$. With the
multiplicity spaces $M_{\pi,\sigma}=\Hom_{G\times H}(V_\pi\otimes W_\sigma,\Omega)$ and
$m_{\pi,\sigma}=\dim M_{\pi,\sigma}$, the $G\times H$-module decomposes as
\begin{equation}\label{eq:product-decomp}
 \Omega
 \cong\bigoplus_{\pi,\sigma}
 V_\pi\otimes W_\sigma\otimes M_{\pi,\sigma},
\end{equation}
and $\Theta(\pi)\cong\bigoplus_\sigma W_\sigma\otimes M_{\pi,\sigma}$ as $H$-modules. The \emph{support} is the
set of pairs with $m_{\pi,\sigma}>0$, and the weighted support graph is as in the introduction. For vector
spaces $A$ and $B$, a vector of $A\otimes B$ is \emph{rank one for the flattening $A\mid B$} if it is a
nonzero decomposable tensor $a\otimes b$. Thus $\mathcal C_\pi^G$ consists of the unit vectors in $\Omega[\pi]$ that
are rank one for $V_\pi\mid\Theta(\pi)$.

The functional $\Phi_K$ is invariant under scalar phase. It is therefore natural to remove this unavoidable
direction before discussing intersections. Put
\[
 X_\pi^G=\mathcal C_\pi^G/S^1\subset\PP(\Omega),
 \qquad
 X_\sigma^H=\mathcal C_\sigma^H/S^1\subset\PP(\Omega).
\]
Here $\mathcal C_\sigma^H$ denotes the $\sigma$-component of
$\operatorname{LocMax}(\Phi_H)$, and each quotient is
identified with the corresponding set of lines in
$\mathbb{P}(\Omega)$.

We use the same notation $\Phi_K$ for the induced function
on projective space, evaluated at any unit representative:
\[
 \Phi_K([v])
 =
 \Phi_K\!\left(\frac{v}{\|v\|}\right),
 \qquad v\neq 0.
\]
This is well defined because $\Phi_K$ is invariant
under scalar phase.

\begin{theorem}[Local-maximal incidence]\label{thm:incidence}
For every pair $(\pi,\sigma)$,
\[
 X_\pi^G\cap X_\sigma^H\ne\varnothing
 \quad\Longleftrightarrow\quad
 M_{\pi,\sigma}\ne0.
\]
If $M_{\pi,\sigma}\ne0$, then
\begin{equation}\label{eq:intersection-product}
 X_\pi^G\cap X_\sigma^H
 \cong
 \PP(V_\pi)\times\PP(W_\sigma)\times\PP(M_{\pi,\sigma}).
\end{equation}
On the corresponding unit vectors,
\begin{equation}\label{eq:heights}
 \Phi_G=\frac1{d_\pi},
 \qquad
 \Phi_H=\frac1{e_\sigma}.
\end{equation}
Hence
\begin{equation}\label{eq:multiplicity-from-dim}
 m_{\pi,\sigma}
 =\dim_{\C}(X_\pi^G\cap X_\sigma^H)+3-d_\pi-e_\sigma.
\end{equation}

Each point of the intersection determines a cyclic
$G\times H$-submodule of $\Omega$.

Moreover, the map
\begin{equation}\label{eq:cyclic-map}
 [v]\longmapsto \Span_{\C}\{(g,h)v:g\in G,\ h\in H\}
\end{equation}
from the intersection to the Grassmannian of $d_\pi e_\sigma$-dimensional subspaces has image naturally
identified with $\PP(M_{\pi,\sigma})$ and fiber
$\PP(V_\pi)\times\PP(W_\sigma)$.
\end{theorem}

\begin{proof}
A point of $X_\pi^G$ has a representative that is rank one for the flattening
\[
 V_\pi\mid\bigoplus_{\sigma'} W_{\sigma'}\otimes M_{\pi,\sigma'},
\]
and similarly for $X_\sigma^H$. A point in both components therefore lies in the joint block
$V_\pi\otimes W_\sigma\otimes M_{\pi,\sigma}$ and is rank one for both flattenings
\[
 V_\pi\mid(W_\sigma\otimes M_{\pi,\sigma}),
 \qquad
 W_\sigma\mid(V_\pi\otimes M_{\pi,\sigma}).
\]

Write $v=x\otimes u$ using the first rank-one condition,
where $x\neq 0$ and $u\in W_\sigma\otimes M_{\pi,\sigma}$.
The second rank-one condition allows us to write
\[
 v=\sum_i a_i\otimes y\otimes z_i
\]
for some nonzero $y\in W_\sigma$, with
$a_i\in V_\pi$ and $z_i\in M_{\pi,\sigma}$.
Choose $f\in V_\pi^*$ such that $f(x)=1$.
Contracting the $V_\pi$-factor gives
\[
 u
 =
 (f\otimes I_{W_\sigma}\otimes I_{M_{\pi,\sigma}})(v)
 =
 y\otimes z,
 \qquad
 z:=\sum_i f(a_i)z_i.
\]
Since $v\neq 0$, we have $z\neq 0$.
Conversely, every nonzero tensor $x\otimes y\otimes z$
satisfies both rank-one conditions.
Thus the common points are precisely the lines of
fully decomposable tensors $x\otimes y\otimes z$, that is, the Segre variety of the triple tensor product,
the image of $\PP(V_\pi)\times\PP(W_\sigma)\times\PP(M_{\pi,\sigma})$ under
$([x],[y],[z])\mapsto[x\otimes y\otimes z]$ \cite{Landsberg}, which gives \eqref{eq:intersection-product}.
The height formula is the one-group local-maximal value from \cite[Theorem~3]{HungKang2026}, and
\eqref{eq:multiplicity-from-dim} follows by taking complex dimensions.

Finally, irreducibility makes $x$ and $y$ cyclic, so
\[
 \Span_{\C}\{(g,h)(x\otimes y\otimes z):g\in G,h\in H\}
 =V_\pi\otimes W_\sigma\otimes\C z.
\]

The generated submodule depends only on $[z]$,
and distinct lines in $M_{\pi,\sigma}$ give distinct submodules.
Hence the image is naturally identified with
$\mathbb{P}(M_{\pi,\sigma})$.
Under this identification and the Segre identification
of the domain, the map becomes the projection
\[
 ([x],[y],[z])\longmapsto[z].
\]
Its fiber over $[z]$ is therefore
$\mathbb{P}(V_\pi)\times\mathbb{P}(W_\sigma)$,
proving the last assertion.
\end{proof}

The theorem turns the two projective local-maximal decompositions into a weighted bipartite graph: its left
vertices are the $X_\pi^G$, its right vertices are the $X_\sigma^H$, and an edge is present exactly when the
corresponding components meet. The edge weight is $m_{\pi,\sigma}$. Thus this variational incidence graph
is exactly the weighted support graph of \eqref{eq:product-decomp}.

\begin{example}[Four lines and four intersection points]\label{ex:square}
Let $G=H=C_2$, with nonidentity elements $g$ and $h$, respectively.
Let $\Omega$ have orthonormal basis $e_{ij}$, $i,j\in\{0,1\}$, with
\[
 g e_{ij}=(-1)^i e_{ij},
 \qquad
 h e_{ij}=(-1)^j e_{ij}.
\]
The two projective $G$-local-maximal components and the two projective $H$-components are lines:
\[
 X_i^G=\PP\Span\{e_{i0},e_{i1}\},
 \qquad
 X_j^H=\PP\Span\{e_{0j},e_{1j}\}.
\]
They meet in four points,
\[
 X_i^G\cap X_j^H=\{[e_{ij}]\},
\]
so the incidence graph is the square $K_{2,2}$. If
$v=ae_{00}+be_{01}+ce_{10}+de_{11}$ is unit, then
\[
 \Phi_G(v)
 =(|a|^2+|b|^2)^2+(|c|^2+|d|^2)^2,
\]
and the analogous column formula holds for $\Phi_H$. Thus restricting one functional to a local-maximal
line of the other picks out exactly its two intersection points. By Theorem~\ref{thm:positive-sum}, the
same four points are the projective local maxima of $s\Phi_G+t\Phi_H$ for every $s,t>0$.
\end{example}

\begin{remark}[Residual symmetry]\label{rem:residual}
On one joint block, Schur's lemma gives
\[
 \End_{G\times H}
 \bigl(V_\pi\otimes W_\sigma\otimes M_{\pi,\sigma}\bigr)
 =I\otimes I\otimes\End(M_{\pi,\sigma}).
\]
The unitary group of the last factor acts transitively on $\PP(M_{\pi,\sigma})$ while preserving both
variational landscapes. Thus multiplicity naturally appears as a projective family of copies; when
$m_{\pi,\sigma}>1$, the two functionals alone cannot distinguish a preferred copy.
\end{remark}

The incidence theorem makes the usual Howe condition particularly transparent.

\begin{corollary}[Howe matching]\label{cor:howe}
Let $\mathcal A_G,\mathcal A_H\subset\End(\Omega)$ be the finite-dimensional algebras generated by
the two group actions. They are mutual commutants if and only if the weighted incidence graph is a matching and
every edge has weight one. Equivalently,
\[
 \Omega=\bigoplus_\alpha V_\alpha\otimes W_\alpha
\]
where the $V_\alpha$ are pairwise inequivalent irreducible
$G$-modules and the $W_\alpha$ are pairwise inequivalent
irreducible $H$-modules.
On a paired block, if
\[
 v=\sum_j s_jx_j\otimes y_j,
 \qquad \sum_j s_j^2=1,
\]
is a Schmidt decomposition, then
\begin{equation}\label{eq:schmidt-purity}
 (\dim V_\alpha)\Phi_G(v)
 =(\dim W_\alpha)\Phi_H(v)
 =\sum_j s_j^4,
\end{equation}
and
\begin{equation}\label{eq:schmidt-map}
 T_G(v)=T_H(v)
 =\frac{\sum_j s_j^3x_j\otimes y_j}{(\sum_j s_j^6)^{1/2}}.
\end{equation}
In particular,
\[
 \LocMax(\Phi_G)=\LocMax(\Phi_H).
\]
\end{corollary}

\begin{proof}
By Burnside's theorem,
$\mathcal A_G=\bigoplus_\pi\End(V_\pi)\otimes I$, so
$\mathcal A_G'=\bigoplus_\pi I\otimes\End(\Theta(\pi))$ \cite{GoodmanWallach}. The subalgebra
$\mathcal A_H$ equals this commutant exactly when each $\Theta(\pi)$ is an irreducible $H$-module and the
modules $\Theta(\pi)$ are pairwise inequivalent, that is, when every $G$-type has one $H$-partner, every
$H$-type has one $G$-partner, and every product multiplicity is one. This is precisely the stated weighted
matching condition.

Under this matching condition, Schur's lemma and the
pairwise inequivalence of the irreducible $H$-modules
$W_\alpha$ give
\[
 \mathcal{A}_H'
 =
 \bigoplus_\alpha
 \operatorname{End}(V_\alpha)\otimes I_{W_\alpha}
 =
 \mathcal{A}_G.
\]
Together with $\mathcal{A}_H=\mathcal{A}_G'$,
this proves that the two algebras are mutual commutants.

On a paired block, Schur averaging gives
\[
 \mathcal R_G(P_v)
 =\frac{I_{V_\alpha}}{\dim V_\alpha}\otimes\sum_j s_j^2P_{y_j},
 \qquad
 \mathcal R_H(P_v)
 =\left(\sum_j s_j^2P_{x_j}\right)\otimes\frac{I_{W_\alpha}}{\dim W_\alpha}.
\]
Taking squared Hilbert--Schmidt norms proves \eqref{eq:schmidt-purity}; applying the two operators to $v$ and
normalizing proves \eqref{eq:schmidt-map}.

By \eqref{eq:criterion}, the local maxima of either
functional on $S(\Omega)$ are precisely the unit
decomposable tensors in the individual paired blocks
$V_\alpha\otimes W_\alpha$.
Hence
\[
 \operatorname{LocMax}(\Phi_G)
 =
 \operatorname{LocMax}(\Phi_H).
\]
\end{proof}

Equality of the two local-maximal sets is slightly weaker than Howe matching, and the incidence
geometry shows exactly by how much.

\begin{proposition}[When the local-maximal loci coincide]\label{prop:equal-locmax}
One has $\LocMax(\Phi_G)=\LocMax(\Phi_H)$ if and only if the support is a matching and
$m_{\pi,\sigma}=1$ for every edge with $d_\pi e_\sigma>1$. In particular, if every constituent
$V_\pi\otimes W_\sigma$ has dimension greater than one, equality of the two local-maximal sets is
equivalent to mutual commutants.
\end{proposition}

\begin{proof}
Assume first that $\LocMax(\Phi_G)=\LocMax(\Phi_H)$. Then every point of every
$X_\pi^G$ lies in some $H$-local-maximal component. If $\pi$ has two partners
$\sigma_1,\sigma_2$, then
$X_\pi^G$ contains $[x\otimes(y_1\otimes z_1+y_2\otimes z_2)]$ with both summands nonzero, which lies in
no $H$-component; symmetrically on the other side. Hence equality forces a matching, and then it
amounts to $X_\pi^G=X_\sigma^H$ for each edge, that is, to equality of the Segre varieties
$\PP(V_\pi)\times\PP(W_\sigma\otimes M_{\pi,\sigma})$ and
$\PP(W_\sigma)\times\PP(V_\pi\otimes M_{\pi,\sigma})$ inside
$\PP(V_\pi\otimes W_\sigma\otimes M_{\pi,\sigma})$. If $m_{\pi,\sigma}=1$, both equal
$\PP(V_\pi)\times\PP(W_\sigma)$. If $d_\pi=e_\sigma=1$, both equal $\PP(M_{\pi,\sigma})$. If, say,
$d_\pi=1<e_\sigma$ and $m_{\pi,\sigma}\ge2$, the first is all of $\PP(W_\sigma\otimes M_{\pi,\sigma})$
while the second is a proper Segre subvariety; the case $d_\pi,e_\sigma\ge2$, $m_{\pi,\sigma}\ge2$ is
similar, since then neither variety contains the other. Conversely, if the support is a matching and
$m_{\pi,\sigma}=1$ whenever $d_\pi e_\sigma>1$, then on every edge either $m_{\pi,\sigma}=1$ or
$d_\pi=e_\sigma=1$, and the two local-maximal loci on that block are equal.
Hence the two local-maximal loci coincide.

The final assertion follows from the preceding corollary,
since the additional dimension hypothesis excludes
the one-dimensional exception.
\end{proof}

\section{Compatible positive sums}\label{sec:positive-sum}

Compatible local-maximal conditions can themselves be realized as the local maxima of one functional.

\begin{theorem}[Positive-sum principle]\label{thm:positive-sum}
Let $K_1,\ldots,K_s$ be finite or compact groups acting unitarily on a nonzero finite-dimensional
Hermitian space $E$.
Assume either that the actions commute pairwise, or that
the groups form a chain $K_1\leq\cdots\leq K_s$
and their actions are restrictions of a single unitary
representation of $K_s$ on $E$.
For arbitrary $a_i>0$, put
\[
 \Psi_a=\sum_{i=1}^s a_i\Phi_{K_i}.
\]
Then
\begin{equation}\label{eq:positive-sum-locmax}
 \LocMax(\Psi_a)=\bigcap_{i=1}^s\LocMax(\Phi_{K_i}).
\end{equation}
At every critical point of $\Psi_a$ that is not a local maximum, its Hessian on the unit sphere
(with respect to the round metric) has a strictly positive direction.
\end{theorem}

\begin{proof}
Only the inclusion $\LocMax(\Psi_a)\subseteq\bigcap_i\LocMax(\Phi_{K_i})$ requires proof.

Suppose $v\in S(E)$ is not a local maximum of $\Phi_{K_j}$
for at least one index $j$.
For commuting actions, choose any such $j$ and take $U=E$.
For a subgroup chain, choose $j$ maximal with this property.
If $j<s$, set $U=\operatorname{span}(K_{j+1}v)$;
by the maximality of $j$ and \eqref{eq:criterion},
this is an irreducible $K_{j+1}$-module.
If $j=s$, take $U=E$.

Using a Schmidt decomposition in each $K_j$-isotypic
component of $U$, write
\[
 v=\sum_\ell\sqrt{p_\ell}\,\xi_\ell,
 \qquad
 \xi_\ell=x_\ell\otimes n_\ell,
\]
where the sum runs over the nonzero Schmidt terms.
The vectors $\xi_\ell$ are orthonormal,
$p_\ell>0$, and $\sum_\ell p_\ell=1$.
Since $\operatorname{span}(K_jv)$ is reducible,
at least two masses are nonzero.
Choose distinct indices $r,t$ and define
\[
 v(\tau)
 =
 \sqrt{p_r+\tau}\,\xi_r
 +
 \sqrt{p_t-\tau}\,\xi_t
 +
 \sum_{\ell\neq r,t}\sqrt{p_\ell}\,\xi_\ell,
 \qquad
 -p_r<\tau<p_t.
\]
This is a smooth curve in $S(U)$ with $v(0)=v$.
Along this curve,
$p_r(\tau)=p_r+\tau$ and $p_t(\tau)=p_t-\tau$,
while all singular directions and the remaining masses
are kept fixed.

If the two masses belong to $K_j$-types of degrees $d_{\alpha_r}$ and $d_{\alpha_t}$, then
\eqref{eq:singular-formula} gives
\begin{equation}\label{eq:strict-mass-transfer}
 \frac{d^2}{d\tau^2}\Phi_{K_j}(v(\tau))
 =2\left(\frac1{d_{\alpha_r}}+\frac1{d_{\alpha_t}}\right)>0.
\end{equation}

Consider $i\neq j$ in the commuting case and $i<j$
in the chain case.
In either setting, $K_i$ preserves the $K_j$-isotypic
components of $U$.
Thus, for $k\in K_i$, mixed terms between distinct
$K_j$-isotypic components vanish.

Within a fixed $K_j$-isotypic component of type $\alpha$,
the Schmidt vectors
$\xi_\ell=x_\ell\otimes n_\ell$ have both
$\{x_\ell\}$ and $\{n_\ell\}$ orthonormal.
In the commuting case, $k$ acts on this block as
$I\otimes\eta_{i,\alpha}(k)$ for a representation
$\eta_{i,\alpha}$ on the multiplicity space.
Hence, for $\ell\neq q$,
\[
 \langle\xi_\ell,k\xi_q\rangle
 =
 \langle x_\ell,x_q\rangle
 \langle n_\ell,\eta_{i,\alpha}(k)n_q\rangle
 =
 0.
\]
In the chain case, $k$ acts as $\alpha(k)\otimes I$,
so instead
\[
 \langle\xi_\ell,k\xi_q\rangle
 =
 \langle x_\ell,\alpha(k)x_q\rangle
 \langle n_\ell,n_q\rangle
 =
 0.
\]
Consequently, in either case,
\[
 \begin{aligned}
 \langle v(\tau),kv(\tau)\rangle
 &=
 \sum_\ell p_\ell(\tau)
 \langle\xi_\ell,k\xi_\ell\rangle\\
 &=
 A_i(k)+\tau D_i(k),
 \end{aligned}
\]
where $A_i(k)=\langle v,kv\rangle$ and
\[
 D_i(k)
 =
 \langle\xi_r,k\xi_r\rangle
 -
 \langle\xi_t,k\xi_t\rangle.
\]
It follows that
\[
 \frac{d^2}{d\tau^2}\Phi_{K_i}(v(\tau))
 =
 2\int_{K_i}|D_i(k)|^2\,d\mu_i(k)
 \geq 0,
\]
where $\mu_i$ is normalized Haar measure,
or normalized counting measure when $K_i$ is finite.

For a chain and $i>j$, the curve stays inside the irreducible $K_i$-module $\Span(K_i v)$, so
$\Phi_{K_i}$ is constant. Therefore $\Psi_a(v(\tau))$ is strictly convex, and $v$ cannot be a local
maximum. Conversely, a common local maximum clearly maximizes every positive sum locally, proving
\eqref{eq:positive-sum-locmax}. Finally, the curve $v(\tau)$ is smooth at $\tau=0$ with
$v'(0)=\xi_r/(2\sqrt{p_r})-\xi_t/(2\sqrt{p_t})\ne0$. At a critical point of $\Psi_a$ the Riemannian
gradient vanishes, so the second derivative of $\Psi_a(v(\tau))$ at $\tau=0$ equals the Hessian of
$\Psi_a$ on the sphere evaluated at $v'(0)$; it is therefore a strictly positive Hessian direction.
\end{proof}

\begin{example}[Necessity of compatibility]
\label{ex:necessity-of-compatibility}
The compatibility hypothesis is essential.
Consider the two order-two actions on $\mathbb{C}^2$
generated by
\[
 Z=
 \begin{pmatrix}
 1 & 0\\
 0 & -1
 \end{pmatrix},
 \qquad
 X=
 \begin{pmatrix}
 0 & 1\\
 1 & 0
 \end{pmatrix}.
\]
Their projective local-maximal loci are the respective
eigenline sets
\[
 \{[1:0],[0:1]\}
 \qquad\text{and}\qquad
 \{[1:1],[1:-1]\},
\]
which are disjoint.
For equal weights and a unit vector $v=(z_1,z_2)$, however,
\[
 \Phi_{\langle Z\rangle}(v)
 +
 \Phi_{\langle X\rangle}(v)
 =
 \frac{3}{2}
 -
 2\bigl(\operatorname{Im}(\overline{z_1}z_2)\bigr)^2.
\]
Thus the sum has an entire circle of projective maxima,
whereas the two summands have no common local maximum.
Hence the positive-sum identity can fail without
the compatibility hypothesis.
\end{example}

\section{Iterated incidence along a branching chain}\label{sec:bratteli}

The two-level subgroup incidence formula of \cite[Theorem~6]{HungKang2026} can be iterated
along a chain. This gives a variational realization of the weighted branching diagram defined below,
including branching multiplicities.

Let
\[
 G_0\leq G_1\leq\cdots\leq G_n
\]
be finite or compact groups.
Fix an irreducible type
$\lambda_n\in\operatorname{Irr}(G_n)$,
and let $V_{\lambda_n}$ be a representative
irreducible unitary $G_n$-module.
For $1\leq j\leq n$ and irreducible types
$\lambda_j\in\operatorname{Irr}(G_j)$ and
$\lambda_{j-1}\in\operatorname{Irr}(G_{j-1})$,
write $V_{\lambda_j}$ and $V_{\lambda_{j-1}}$
for representative irreducible unitary modules, and put
\begin{equation}\label{eq:branching-space}
 B(\lambda_j,\lambda_{j-1})
 :=
 \Hom_{G_{j-1}}
 \bigl(V_{\lambda_{j-1}},
       \Res^{G_j}_{G_{j-1}}V_{\lambda_j}\bigr).
\end{equation}
A \emph{Bratteli path} below $\lambda_n$ is a sequence
\[
 p=(\lambda_n,\lambda_{n-1},\ldots,\lambda_0)
\]
for which every space in \eqref{eq:branching-space} is nonzero. Equivalently, it is a downward path from
$\lambda_n$ in the \emph{weighted Bratteli diagram}: the graded graph whose vertices at level $j$ are the
irreducible representations of $G_j$, with an edge from $\lambda_j$ to $\lambda_{j-1}$ of weight
$\dim B(\lambda_j,\lambda_{j-1})$ whenever this weight is positive. We keep multiplicity as an edge
weight rather than replacing a weight-$m$ edge by $m$ parallel edges; thus a path here means a vertex path
in this weighted diagram. Let $X_p\subset\PP(V_{\lambda_n})$
be the set of lines whose unit representatives lie in the $\lambda_j$-component of the
local-maximal set for $G_j$ at every level $j$.

\begin{theorem}[Variational Bratteli theorem]\label{thm:bratteli}
For every Bratteli path
$p=(\lambda_n,\ldots,\lambda_0)$, the map
\begin{equation}\label{eq:bratteli-map}
 \PP(V_{\lambda_0})
 \times
 \prod_{j=1}^n\PP\bigl(B(\lambda_j,\lambda_{j-1})\bigr)
 \longrightarrow
 \PP(V_{\lambda_n}),
 \qquad
 ([x],[\varphi_1],\ldots,[\varphi_n])
 \longmapsto
 [\varphi_n\cdots\varphi_1(x)]
\end{equation}
is an embedding with image $X_p$. The simultaneous projective local-maximal locus is the disjoint
union of these $X_p$. On $X_p$,
\begin{equation}\label{eq:bratteli-heights}
 \Phi_{G_j}=\frac1{\dim V_{\lambda_j}},
 \qquad 0\le j\le n.
\end{equation}
\end{theorem}

\begin{proof}
Restrict $V_{\lambda_n}$ to $G_{n-1}$:
\[
 V_{\lambda_n}
 \cong
 \bigoplus_{\lambda_{n-1}}
 V_{\lambda_{n-1}}\otimes
 B(\lambda_n,\lambda_{n-1}).
\]
By the one-group criterion, a $G_{n-1}$-local maximum in the
$\lambda_{n-1}$-isotypic component is a simple tensor
$x_{n-1}\otimes b_n$. Equivalently, it selects an embedded irreducible copy
$\varphi_n(V_{\lambda_{n-1}})$, where
$[\varphi_n]\in\PP(B(\lambda_n,\lambda_{n-1}))$, together with the line
$[\varphi_n(x_{n-1})]$ inside that copy.

The group $G_{n-2}$ acts only on the factor $V_{\lambda_{n-1}}$. Applying the same argument to
$\Res^{G_{n-1}}_{G_{n-2}}V_{\lambda_{n-1}}$ selects a type $\lambda_{n-2}$ and a line
$[\varphi_{n-1}]\in\PP(B(\lambda_{n-1},\lambda_{n-2}))$. Repeating down the chain gives exactly
\eqref{eq:bratteli-map}.

Conversely, let $v=\varphi_n\cdots\varphi_1(x)$,
where $x\neq 0$ and each $\varphi_j$ is nonzero.
Each $\varphi_j$ is injective, since its domain is irreducible.
For every $j$, the map $\varphi_n\cdots\varphi_{j+1}$
is $G_j$-equivariant, and the nonzero vector
$\varphi_j\cdots\varphi_1(x)$ generates $V_{\lambda_j}$.
Consequently,
\[
 \operatorname{span}(G_jv)
 =
 (\varphi_n\cdots\varphi_{j+1})(V_{\lambda_j})
 \cong V_{\lambda_j},
 \qquad 0\leq j\leq n,
\]
where empty compositions are interpreted as identity maps.
By \eqref{eq:criterion}, the unit normalization of $v$
is therefore locally maximal for every $\Phi_{G_j}$.

The projective point $[v]$ determines the cyclic submodule
\[
 \operatorname{span}(G_{n-1}v)
 =
 \varphi_n(V_{\lambda_{n-1}}).
\]
By Schur's lemma, this submodule determines the
intertwiner line $[\varphi_n]$ uniquely.
Pulling $[v]$ back along $\varphi_n$ gives a projective
point in $V_{\lambda_{n-1}}$, independently of the
choice of representative of $[\varphi_n]$.
Repeating this procedure down the chain recovers
all the remaining intertwiner lines $[\varphi_j]$
and finally the bottom line $[x]$.
Thus the map \eqref{eq:bratteli-map} is injective.

Iterating the canonical evaluation isomorphisms gives a linear isomorphism from the direct sum over vertex
paths onto $V_{\lambda_n}$; on the $p$-summand it is exactly
$V_{\lambda_0}\otimes B(\lambda_1,\lambda_0)\otimes\cdots\otimes B(\lambda_n,\lambda_{n-1})
\longrightarrow V_{\lambda_n}$,
$x\otimes\varphi_1\otimes\cdots\otimes\varphi_n\mapsto\varphi_n\cdots\varphi_1(x)$.
Hence \eqref{eq:bratteli-map} is the composition of a Segre embedding with the projectivization of an
injective linear map, and is therefore an embedding of projective varieties. The selected isomorphism type at each level makes different path components
disjoint. Formula \eqref{eq:bratteli-heights} is
the one-group local-maximal value. Each $X_p$ is connected and closed, and only finitely many
paths occur below $\lambda_n$; hence these are exactly the connected components of the simultaneous locus.
\end{proof}

Combining Theorems~\ref{thm:positive-sum} and \ref{thm:bratteli}, for arbitrary weights
$a_0,\ldots,a_n>0$ the projective local-maximal components of
\[
 \Psi_a=\sum_{j=0}^n a_j\Phi_{G_j}
\]
are exactly the varieties $X_p$ indexed by Bratteli paths. On the component for
$p=(\lambda_n,\ldots,\lambda_0)$ the value is
\begin{equation}\label{eq:path-sum-height}
 \Psi_a|_{X_p}=\sum_{j=0}^n\frac{a_j}{\dim V_{\lambda_j}}.
\end{equation}
Thus the whole weighted branching geometry is realized by one scalar quartic functional.

A basis of $V_{\lambda_n}$ is \emph{adapted to the chain}
if each basis vector generates an irreducible
$G_j$-submodule for every $0\leq j\leq n$.

\begin{corollary}[Variational basis criterion]\label{cor:bratteli-mf}
For a Bratteli path $p=(\lambda_n,\ldots,\lambda_0)$,
the component $X_p$ is a point if and only if
\[
 \dim V_{\lambda_0}=1
 \qquad\text{and}\qquad
 \dim B(\lambda_j,\lambda_{j-1})=1
 \quad (1\leq j\leq n).
\]
If these conditions hold for every path below $\lambda_n$,
write $X_p=\{[v_p]\}$.
Choosing each $v_p$ to have unit norm gives an
orthonormal basis of $V_{\lambda_n}$ adapted to the chain.
\end{corollary}

\begin{proof}
The first statement follows immediately from the product description in
Theorem~\ref{thm:bratteli}. Distinct paths first separate inside orthogonal irreducible summands at some
level, so their lines are orthogonal. The branching dimension recursion, with one-dimensional bottom
types and edges, shows that the number of paths is $\dim V_{\lambda_n}$.
\end{proof}

The resulting basis is the Gelfand--Tsetlin basis
of the chain, unique up to independent phases.

\begin{remark}[Residual multiplicity along a path]\label{rem:path-multiplicity}
Without the multiplicity-free hypothesis, Theorem~\ref{thm:bratteli} identifies exactly what remains:
each branching edge of multiplicity $m$ contributes a factor $\PP^{m-1}$, while the bottom type contributes
$\PP(V_{\lambda_0})$. Thus a path reduces to a distinguished line precisely when every branching
multiplicity and the bottom dimension are one.
\end{remark}

\section{Classical dualities}\label{sec:examples}

We now specialize Theorems~\ref{thm:incidence} and \ref{thm:bratteli} to three classical pictures:
Peter--Weyl and Frobenius reciprocity, Gelfand--Tsetlin bases, and Schur--Weyl duality. The first reads
pairwise incidence as restriction multiplicity; the last two show how a subgroup chain resolves residual
multiplicity into the familiar adapted bases.

\subsection{Regular representations and restriction}

Let $G$ be a finite group and $K\leq G$.
Identify $\Omega=\mathbb{C}[G]$ with the space of
complex-valued functions on $G$, equipped with the
normalized counting inner product.
Write $G_L$ and $G_R$ for two copies of $G$ acting by
left and right translation:
\[
 (L_gf)(x)=f(g^{-1}x),
 \qquad
 (R_hf)(x)=f(xh),
 \qquad g,h,x\in G.
\]
These actions commute.
Write $K_R\leq G_R$ for the copy of $K$ acting by
restriction of the right translation action.
For $\tau\in\operatorname{Irr}(K)$, let $W_\tau$
be its representation space.

Recall the finite-group Peter--Weyl decomposition \cite{GoodmanWallach}
\[
 \mathbb{C}[G]
 \cong
 \bigoplus_{\pi\in\operatorname{Irr}(G)}
 V_\pi\otimes V_\pi^*,
\]
where $G_L$ acts on the first factor and $G_R$
acts on the second factor by the dual representation.
The $\pi$-block is spanned by the functions
\[
 c_{u,w}(x)=\langle u,\pi(x^{-1})w\rangle,
 \qquad u,w\in V_\pi.
\]

\begin{example}[Peter--Weyl]\label{ex:PW}
For $K=G$ the support is the matching $\pi\leftrightarrow\pi^*$ with all weights one, so by
Corollary~\ref{cor:howe} the two local-maximal loci coincide. Their common projective
local-maximal component is
\[
 X_\pi^{G_L}=X_{\pi^*}^{G_R}
 \cong\PP(V_\pi)\times\PP(V_\pi^*),
\]
the Segre variety of rank-one elements of $\End(V_\pi)$, at height $1/d_\pi$. In words: \emph{the local
maxima of the matrix-coefficient functional on the regular representation are precisely
the unit normalizations of nonzero single
matrix-coefficient functions $c_{u,w}$.} For $G=S_3$ and its two-dimensional representation this is the quadric
$\PP^1\times\PP^1\subset\PP^3$.
\end{example}

\begin{example}[Branching and Frobenius reciprocity]\label{ex:branching}
Recall that Frobenius reciprocity \cite{GoodmanWallach} gives
\[
 \operatorname{Hom}_G
 \bigl(U,\operatorname{Ind}_K^G W\bigr)
 \cong
 \operatorname{Hom}_K
 \bigl(\operatorname{Res}_K^G U,W\bigr)
\]
for a $G$-module $U$ and a $K$-module $W$.
Here $\operatorname{Ind}_K^G W$ denotes the representation
of $G$ induced from $W$.
Thus induction and restriction encode the same
multiplicity data.

For a proper subgroup $K$, the $(\pi,\tau)$-multiplicity space is $\Hom_K(W_\tau,V_\pi^*)$, so
\[
 X_\pi^{G_L}\cap X_\tau^{K_R}
 \cong
 \PP(V_\pi)\times\PP(W_\tau)\times
 \PP\bigl(\Hom_K(W_\tau,V_\pi^*)\bigr).
\]

Thus the incidence graph is the branching graph for
restriction from $G$ to $K$, with the left vertex $\pi$
corresponding to the contragredient representation $\pi^*$.
The branching multiplicity is recovered from the dimension
of the intersection by Theorem~\ref{thm:incidence}.

Read from the $K_R$-side, the residual $G_L$-module
of Proposition~\ref{prop:nested} is
\[
 \operatorname{Hom}_{K_R}
 \bigl(W_\tau,\mathbb{C}[G]\bigr)
 \cong
 \operatorname{Ind}_K^G W_\tau^*,
\]
using the usual function model for induction.
Frobenius reciprocity and duality then give
\[
 \begin{aligned}
 \operatorname{Hom}_G
 \bigl(V_\pi,\operatorname{Ind}_K^G W_\tau^*\bigr)
 &\cong
 \operatorname{Hom}_K
 \bigl(\operatorname{Res}_K^G V_\pi,W_\tau^*\bigr)\\
 &\cong
 \operatorname{Hom}_K
 \bigl(W_\tau,\operatorname{Res}_K^G V_\pi^*\bigr).
 \end{aligned}
\]
Thus the same multiplicity space appears both in the restriction of $V_\pi^*$ to $K$ and in
the decomposition of the residual $G_L$-module. This gives a left--right regular-representation
realization of the branching--induction correspondence discussed in
\cite[Section~3.2]{HungKang2026}.

For $S_n\supset S_{n-1}$, whose representations are self-dual and whose branching
is multiplicity-free, every nonempty intersection is a plain Segre product $\PP(V_\lambda)\times\PP(V_\mu)$,
and the incidence graph is the part of Young's graph between levels $n$ and $n-1$. Here \emph{Young's
graph} has the partitions as vertices, with $\lambda\to\mu$ when $\mu$ is obtained from $\lambda$ by removing
one box; it is the Bratteli diagram of the chain $S_1\subset S_2\subset\cdots$.
\end{example}

\subsection{Gelfand--Tsetlin bases and residual multiplicity}

For the symmetric-group chain
\[
 S_1\subset S_2\subset\cdots\subset S_n,
\]
Young branching is multiplicity-free and the bottom $S_1$-type is one-dimensional
\cite{OkounkovVershik}. Hence
Corollary~\ref{cor:bratteli-mf} says that the simultaneous local-maximal lines in an irreducible
$S_n$-module $V_\lambda$ are indexed by paths in Young's graph, equivalently by standard Young tableaux
of shape $\lambda$. Thus Young's orthogonal basis is characterized, up to phases, by the common
projective local maxima of the subgroup functionals. The same argument applies to
$U(1)\subset\cdots\subset U(n)$, since successive restrictions are multiplicity-free and all irreducible
$U(1)$-modules are one-dimensional \cite{GoodmanWallach}. It gives the lines of the usual
Gelfand--Tsetlin basis. This variational
characterization is different from the classical combinatorial and functional realizations of the
Gelfand--Tsetlin basis; see, for example, \cite{HershLenart,Artamonov}.

\begin{example}[A chain resolves a multiplicity space]\label{ex:S5}
Let $V=V_{(3,2)}$ for $S_5$. In the restriction to $S_3$, the type $(2,1)$ occurs with multiplicity
two, corresponding to the two paths
\[
 (3,2)\to(3,1)\to(2,1),
 \qquad
 (3,2)\to(2,2)\to(2,1).
\]
Hence the projective $S_3$-local-maximal component for this type is
\[
 \PP(V_{(2,1)})\times\PP^1.
\]
The second factor parametrizes the copies of $V_{(2,1)}$ inside $V$. After inserting the intermediate
group $S_4$, Theorem~\ref{thm:bratteli} replaces this continuous copy parameter by the two distinguished
path points:
\[
 \PP(V_{(2,1)})\times\{p_{(3,1)},p_{(2,2)}\}.
\]
Thus the subgroup chain resolves the residual projective multiplicity space exactly according to the
Bratteli paths. Theorem~\ref{thm:positive-sum} makes this resolution visible in one function. For
\[
 \Psi=\Phi_{S_3}+\Phi_{S_4},
\]
the projective local-maximal components are the common projective local-maximal components of
$\Phi_{S_3}$ and $\Phi_{S_4}$, indexed by the partial paths $(3,2)\to\lambda_4\to\lambda_3$.
On these components, $\Psi$ takes the values
\[
 \frac43\quad\text{on }(3,1)\to(3),
 \qquad
 \frac56\quad\text{on }(3,1)\to(2,1),
 \qquad
 1\quad\text{on }(2,2)\to(2,1).
\]
Thus the representation-theoretically meaningful local maxima need not have the same height; only the
first component is globally maximal for this choice of weights.
\end{example}

\subsection{Schur--Weyl duality}\label{sec:SW}

Let $S_k$ permute the factors of $(\C^n)^{\otimes k}$ and let $U(n)$ act diagonally. Schur--Weyl duality
\cite{GoodmanWallach} gives the matching
\[
 (\C^n)^{\otimes k}
 \cong
 \bigoplus_{\lambda\vdash k,\ \ell(\lambda)\le n}
 S_\lambda(\C^n)\otimes S^\lambda,
\]
with $S_\lambda(\C^n)$ the Weyl module, $S^\lambda$ the Specht module, and $\ell(\lambda)$ the number of
parts of $\lambda$. By
Remark~\ref{rem:compact} and Corollary~\ref{cor:howe}, the two local-maximal loci coincide, and the
common projective $\lambda$-component is $\PP(S_\lambda(\C^n))\times\PP(S^\lambda)$.
In particular, the embedded copies of
the Weyl module $S_\lambda(\C^n)$ in $(\C^n)^{\otimes k}$ form the projective space
$\PP(S^\lambda)$. Applying Theorem~\ref{thm:bratteli} to the irreducible $S_k$-module
$S^\lambda$ along the chain $S_1\subset\cdots\subset S_k$ singles out its Young
(Gelfand--Tsetlin) basis lines, indexed by standard tableaux. These lines determine the copies
$S_\lambda(\C^n)\otimes\C b$ used in the quantum Schur transform \cite{BaconChuangHarrow}.

\begin{example}[Three qubits]\label{ex:qubits}
For $n=2$ and $k=3$, the Schur--Weyl decomposition is
$S_{(3)}(\C^2)\otimes S^{(3)}\oplus S_{(2,1)}(\C^2)\otimes S^{(2,1)}$.
In physical language, after restriction to $SU(2)$,
$S_{(3)}(\C^2)=\mathrm{Sym}^3\C^2$ is the spin-$3/2$ representation of $SU(2)$, and
$S_{(2,1)}(\C^2)\cong\C^2$ is spin $1/2$, occurring with the two-dimensional multiplicity space
$S^{(2,1)}$; the copies of spin $1/2$ are therefore parametrized by $\PP(S^{(2,1)})\cong\PP^1$. Since
$S^{(2,1)}$ restricts to $S_2$ as the sum of the trivial and sign characters, Theorem~\ref{thm:bratteli}
selects two points $p_0,p_1$ on this line: the copies on which the transposition of the first two factors
acts by $-1$ or by $+1$, that is, on which the first two spins couple to total spin $0$ or $1$.
Thus the simultaneous projective local-maximal locus for $U(2)$, $S_3$, and $S_2$ in the
spin-$1/2$ sector is
\[
 \PP\bigl(S_{(2,1)}(\C^2)\bigr)\times\{p_0,p_1\}.
\]
\end{example}

\medskip
\noindent
Across these examples, pairwise incidence detects support and multiplicity, while a subgroup chain
resolves the remaining multiplicity into Bratteli paths. The positive-sum principle imposes these
compatible conditions through one scalar functional.

\section*{Declaration of generative AI and AI-assisted technologies in the manuscript preparation process}

During the preparation of this manuscript, the authors used OpenAI ChatGPT and Anthropic Claude to assist with proof checking and language revision. The authors are responsible for reviewing and verifying the mathematical arguments and references, and take full responsibility for the content of the manuscript.

\end{document}